 \documentclass[11pt,activeacute,amsfonts,reqno]{amsart}
\usepackage{todonotes}
\usepackage{appendix}
\usepackage[english]{babel}
\usepackage{inputenc}
\usepackage{mathabx}
\usepackage{amsmath,amsthm,amsfonts,amssymb,amscd}
\usepackage{amssymb,amsfonts}
\usepackage[all,arc]{xy}
\usepackage{enumerate}
\usepackage{mathrsfs}
\usepackage{marginnote}
\usepackage{mathtools}
\usepackage{amsthm,thmtools,xcolor}

\usepackage[all]{xy}
\usepackage{tikz-cd}
\usetikzlibrary{cd}
\usepackage[square,sort,comma,numbers,nonamebreak]{natbib}

\theoremstyle{plain}

\newtheorem{lemx}{Lemma}

\theoremstyle{definition}

\theoremstyle{remark}

\numberwithin{equation}{section}
\usepackage{color}

\usepackage[colorlinks = true,
linkcolor = blue,
urlcolor  = blue,
citecolor = blue,
anchorcolor = blue]{hyperref}
\usepackage{hyperref}
\numberwithin{equation}{section}

\DeclareMathOperator{\sgn}{\mathrm{sgn}}

\newcommand\underrel[3][]{\mathrel{\mathop{#3}\limits_{%
			\ifx c#1\relax\mathclap{#2}\else#2\fi}}}

\usepackage{xcolor}
\usepackage{graphicx}
\usepackage{wrapfig}
\usepackage{tikz}
\usetikzlibrary{arrows,calc,decorations.pathreplacing}
\usepackage{pgfplots}
\usepackage{qtree}
\usepackage{multimedia}
\usepackage[T1]{fontenc}
\usepackage{mathptmx}
\usepackage{hyperref}

\numberwithin{equation}{section}

\newcommand{\R}{\mathbb{R}}

\newcommand{\al}{\alpha}
\newcommand{\bt}{\beta}
\newcommand{\ga}{\gamma}

\def\bm{\left( \begin{array}{cc}}
\def\endm{\end{array}\right)}

\newcommand{\be}{\begin{equation}}
\newcommand{\ee}{\end{equation}}
\newcommand{\ba}{\begin{equation*}}
\newcommand{\ea}{\begin{equation*}}
\newcommand{\bea}{\begin{eqnarray}}
\newcommand{\eea}{\end{eqnarray}}
\newcommand{\bee}{\begin{eqnarray*}}
\newcommand{\eee}{\end{eqnarray*}}
\newcommand{\ben}{\begin{enumerate}}
\newcommand{\een}{\end{enumerate}}

\title[Larga data asymptotics]{Review on long time asymptotics of large data in some nonintegrable dispersive models}

\author[C. Mu\~noz]{Claudio Mu\~noz$^{\lor}$}
\address{CNRS and Departamento de Ingenier\'{\i}a Matem\'atica and Centro
de Modelamiento Matem\'atico (AFB170001 and UMI 2807 CNRS), Universidad de Chile, Casilla
170 Correo 3, Santiago, Chile.}
\email{cmunoz@dim.uchile.cl}
\thanks{$^{\lor}$ C.M.'s work was funded in part by Chilean research grants FONDECYT 1191412, project France-Chile
ECOS-Sud C18E06, MathAmSud EEQUADD II and Centro de Modelamiento Matemático (CMM), ACE210010 and FB210005, BASAL funds for centers of excellence from ANID-Chile.}


\date{\today}

\begin{document}

\begin{abstract}
In this short note we review recent results concerning the long time dynamics of large data solutions to several dispersive models. Starting with the KdV case and ending with the KP models, we review the literature and present new results where virial estimates allow one to prove local energy decay in different regions of space where solitons, lumps or solitary waves are not present. %Particular emphasis will be given to 
%We consider the Kadomtsev-Petviashvili equations posed on $\mathbb{R}^2$. For both models, we provide sequential in time asymptotic descriptions of solutions obtained from arbitrarily large initial data, inside regions of the plane not containing lumps or line solitons, and under minimal regularity assumptions. The proof involves the introduction of two new virial identities adapted to the KP dynamics. This new approach is particularly important in the KP-I case, where no monotonicity property was previously known. The core of our results do not require the use of the integrability of KP and are adaptable to well-posed perturbations.
\end{abstract}

\maketitle
%\tableofcontents

\section{Introduction}

These notes are part of a talk given at the \emph{TYAN Virtual Thematic Workshop in Mathematics}, held on September 27th 2021. Here, I will describe some new results concerning the time decay of solutions to Korteweg-de Vries (KdV) and related equations with rough large data, mainly concentrated in regions where solitons are not present. These are results obtained with my collaborators

\medskip

\begin{itemize}
\item Miguel A. Alejo (U. C\'ordoba, Spain),
\item Fernando Cortez (EPN, Ecuador),
\item Chulkwang Kwak (Ehwa Womans U., South Korea),
\item Argenis Mendez (PUCV Chile),
\item Felipe Poblete (UACh, Chile),
\item Gustavo Ponce (UCSB, USA),
\item Juan C. Pozo (U. Chile, Chile),
\item Jean-Claude Saut (U. Paris-Saclay, France),
\end{itemize}
to whom I deeply acknowledge. 

\medskip

\noindent
{\bf Setting.} First of all, we consider the {\bf generalized KdV equations} posed on $\R:$
\[
\hbox{(gKdV)}\quad 
\begin{cases}
\partial_t u+ \partial_x (\partial_x^2 u + u^p)=0, \\
u=u(t,x)\in\R, \quad t,x\in \R, \quad p=2,3,4.
\end{cases}
\]
Some important remarks are in order:
\ben
\item $p=2,3,4$: equation is globally well-posed in $H^1$ (Kato \cite{Ka}, Kenig, Ponce and Vega \cite{KPV1}).
\medskip
\item $p=2,3$ are {\bf completely integrable models} (KdV and mKdV, respectively), see e.g. \cite{AS} and references therein.
\medskip
\item $p=5$ is the $L^2$ critical case, in the sense that the natural scaling in this case $u(t,x)\mapsto c^{1/2}u(c^3t,cx)$ preserves the $L^2$ norm.
\medskip
\item If $p\geq 5$, blow up may occur (Martel-Merle \cite{MM3}).
\een
Essentially, one has global solutions for data in $H^1$ (or $L^2$) and power $p\leq 4$. For more details, the reader can consult the monograph by Linares and Ponce \cite{LP}. Here, we are interested in the following

\medskip

\noindent
{\bf Key question:} \emph{To describe the asymptotics of large, globally defined solutions in nonintegrable cases (a.k.a. the ``soliton resolution conjecture'') and subcritical nonlinearities.}

\medskip

Somehow this reduces the problem to the case $p=4$ in (gKdV), but we are also interested in proofs that are \emph{independent of the integrable character of the model}. In particular, the results that I shall present here are also valid for some of the integrable models, and are concentrated in the regions where solitons/solitary waves/lumps {\bf are not present}.

\medskip

To try to explain the previous key question, we need some definitions.

\medskip

A {\bf soliton} (or solitary wave)  is a solution to (gKdV) of the form  
\[
u(t,x)=Q_c(x-ct),
\]
with $Q_c\in H^1(\R)$ the unique (up to translations) positive solution to
\[
Q_c'' -c Q_c + Q_c^p=0, \quad c>0.
\]
The soliton resolution conjecture states that, except for some particular cases, globally defined solutions to (gKdV) decompose, as time tends to infinity, as the {\bf sum of solitons and radiation}. By radiation, we broadly mean a linear solution of (gKdV). Mathematically,
\[
u(t,x) \sim_{L^2} \sum_{j\geq 1} Q_{c_j(t)}(x-\rho_j(t)) + \ell(t), \quad t\to +\infty. 
\]
Here, $c_j(t)$ and $\rho_j(t)$ are modulated scaling and shifts, respectively, and $ \ell(t)$ denotes a linear solution, namely $\partial_t \ell + \partial_x^3 \ell =0$ ($p\geq 4$, see explanations below). Several questions appear immediately: are the $c_j(t)$ always different? What about the behavior of $\rho_j(t)$? 

\medskip

These are questions that still remain unsolved if $p=4$, for instance. Indeed, except for integrable cases (to be described below), this conjecture is far from being established in the nonintegrable subcritical setting.

\section{A quick review on the literature} For the sake of completeness, we will consider a general nonlinearity $f(u)$ in (gKdV), which will be specified below. Therefore, we have 
\[
 \hbox{(gKdV)} \qquad \partial_t u+ \partial_x (\partial_x^2 u + f(u))=0.
\]
The following account is by far not 100\% accurate and complete, but it describes in a broad sense the recent literature on the subject.

\medskip

\noindent
{\bf Long time behavior of small solutions:} this corresponds to the case where one assumes small initial data in a well-chosen Banach space. The situation is quite good for powers of the nonlinearity $f(u)\sim u^p$ with $p\geq 3$, but below this range things are quite complicated and essentially only proved via inverse scattering techniques and under weighted norms. 

\medskip

We will say that a solution $u(t)\in H^1$ of (gKdV) scatters to a linear one if there exists $\ell(t)\in H^1$ solution of Airy $\partial_t \ell + \partial_x^3 \ell =0$ such that 
\[
\|u(t)-\ell(t)\|_{H^1} \to 0, \quad t\to +\infty.
\]
This is formally the expected case if $p\geq 4$. Sometimes one needs to perturb $\ell(t)$ to recover the asymptotic behavior at infinity, this is the case of \emph{modified scattering} ($p\leq 3$).

\medskip

Concerning nonintegrable techniques, here one has:
\ben
\item {Kenig-Ponce-Vega} \cite{KPV1}: scattering for small data solutions of the $L^2$ critical gKdV equation ($p=5$). 
\medskip
\item {Ponce-Vega} \cite{PV}: for the case $f(s)=|s|^p$, $p>(9+\sqrt{73})/4\sim 4.39$, small data solutions in $L^{1}\cap H^{2}$ lead to decay, with rate $t^{-1/3}$ (i.e. linear rate of decay). 
\medskip
\item {Christ and Weinstein} \cite{CW}: scattering of small data if $f(s)= |s|^p$, $p>\frac14(23-\sqrt{57}) \sim 3.86$. 
\medskip
\item {Hayashi and Naumkin} \cite{HN1,HN2} studied the case $p>3$, obtaining decay estimates and asymptotic profiles for small data in the classical weighted space $H^{1,1}$. 
\een

\medskip

\noindent
{\bf Long time behavior of small solutions around solitons:} In this case, we assume initial data which is close in some sense to a soliton. Consequently, one seeks here stability and asymptotic stability (AS) of the soliton 
in particular spaces.
\ben
\item Bona-Souganidis-Strauss \cite{BSS}: Stability and instability for the gKdV soliton in the subcritical and supercritical cases. 
\medskip
\item {\color{black} Pego-Weinstein} \cite{PW}: First result of AS gKdV models: asymptotic stability in exponentially weighted spaces.
\medskip
\item {\color{black} Martel-Merle} \cite{MM,MM1,MM2,MMcol}: Asymptotics in the energy space using virial techniques, and studied the {collision problem} in quartic gKdV.
\medskip
\item {\color{black} Tao} \cite{Tao} considered the scattering of data for the quartic KdV in the space $H^1 \cap \dot H^{-1/6}$ around the zero and the soliton solution. 
\medskip
\item The finite energy condition above was then removed by {\color{black} Koch and Marzuola} \cite{KM} by using $U-V$ spaces. 
\medskip
\item {\color{black} C\^ote} \cite{Cote} constructed solutions to the subcritical gKdV equations with a given asymptotical behavior, for $p=4$ (quartic) and $p=5$ ($L^2$ critical).
\medskip
\item {\color{black} Germain, Pusateri and Rousset} \cite{GPR} dealt with the mKdV case ($p=3$) around the zero background and the soliton, by using Fourier techniques and estimates on space-time resonances. See also Harrop-Griffiths \cite{HG}.
\medskip
\item {\color{black} Chen-Liu} \cite{CL}: soliton resolution of mKdV using inverse scattering techniques, in weighted Sobolev spaces. This results also includes the presence of nonzero speed breathers \eqref{breather}, to be explained below.
\een

\medskip

No scattering results seems to hold for the quadratic power ($p=2$), which can be considered as ``supercritical'' in terms of modified scattering.

\section{A First result for the KdV model}

Our first result concerns the integrable quadratic case, which is nothing but KdV. Using integrability techniques, it is possible to give a detailed description of different regions, see e.g. \cite{ES,DVZ}. Assuming that the solution stays in the space $L^\infty(\R;L^1(\R))$, one can prove full decay around zero in the region $|x|\lesssim t^{1/2-}$.

\medskip

{\it
\noindent
{\bf Theorem (G. Ponce-M. 2018 \cite{MuPo1}).} Assume $u\in C(\R;H^1(\R))\cap {\color{black} L^\infty(\R;L^1(\R))}$ solution to (KdV) $p=2$, then  
\[
\lim_{t\rightarrow \infty}\int_{\Omega(t)}u^{2}(x,t)\mathrm{d}x=0, \qquad \Omega(t):=\left\{ x\in\R ~: ~ |x|\lesssim t^{1/2}\log^{-2} t \right\}.
\]
}
There are several comments that are important to address:

\medskip

\noindent
{\bf Remarks.} 
\ben
\item Solitons do satisfy the last hypothesis (being in $L^\infty(\R;L^1(\R))$).
\medskip
\item No size restriction on the data is needed.
\medskip
\item No use of integrable techniques, this result is valid for nonintegrable perturbations of KdV and small $H^1$ data, in the form
\[
u^2 + o_{u\to 0}(u^2).
\]
One of the most important examples is the integrable Gardner model:
\[
\partial_t u + \partial_{x}(\partial_x^2 u +  u^2 + \mu u^3) = 0, \quad  \mu>0,
\]
see \cite{Alejo} for more details.

\medskip

\item This decay result is not valid for $p=3$ (existence of breathers, which are counterexamples, see below) and $p=4$ (scaling problems).
\medskip
\item  {\color{black} Ifrim-Koch-Tataru} \cite{IKT} give a detailed description of asymptotics for weighted small data $O(\varepsilon)$ and times of cubic order $O(\varepsilon^{-3})$.
\een

\medskip

\noindent
{\bf About breathers.} It is well-known that both mKdV and Garner models have stable breather solutions \cite{AS,AM,AM1,AM2,Alejo}, that is to say, \emph{localized in space solutions which are also periodic in time}, up to the symmetries of the equation. An example of these type of solutions is the mKdV breather: for any $\al,\bt>0$, 
\be\label{breather}
\begin{aligned}
B(t,x) := &~{} ~2\sqrt{2} \partial_x \arctan \left(\frac{\beta  \sin(\al (x+\delta t))}{\al \cosh(\beta (x+\ga t))} \right),\\
 \delta= &~{} \al^2 -3\bt ^2, \qquad \gamma= 3\al^2 -\bt ^2,
\end{aligned}
\ee
is a solution of mKdV with nontrivial time-periodic behavior, up to the translation symmetries of the equation. Therefore, generalized KdV equations may have both solitary waves and breathers as well, and both classes of solutions do not decay.

\medskip

\noindent
{\bf About the proof}. Since we want to show that 
\[
\lim_{t\rightarrow \infty}\int_{\Omega(t)}u^{2}(x,t)\mathrm{d}x=0,
\]
the key part of the proof is to show sequential decay of the involved quantity. Consider the functional
\[
\int \psi \left( \frac{x}{\lambda(t)} \right) u(t,x)dx, \quad \psi\sim \tanh, \quad \lambda(t) \sim t^{1/2-},
\]
whose derivative satisfies the lower bound
\[
\frac{d}{dt}\int \psi \left( \frac{x}{\lambda(t)} \right) u(t,x)dx \gtrsim \frac1{\lambda(t)}\int \psi' \left( \frac{x}{\lambda(t)} \right)u^2(t,x)dx - f(t),
\] 
with $f(t)\in L^1(t\gg 1)$. This last virial identity proves that the local in space $L^2$ norm of the solution $u$ has integrability in time properties, in particular, at least for a sequence of times is converging to zero.

\medskip

Part of the problem, consequently, is to remove the additional hypothesis $u\in L^\infty(\R;L^1(\R))$, which is clearly unsatisfactory. Part of the problem was solved for the Benjamin-Ono model, which we describe now.

\section{The Benjamin-Ono model}

This is the model given by
\[
\begin{aligned}
\hbox{(BO)}\quad \begin{cases}
& \!\!  \partial_t u + \partial_{x}(\partial_x\mathcal H u +u^2) = 0,\qquad (t,x) \in \R\times \R,\\
&\!\! u(0,x)=u_0(x),
 \end{cases}
\end{aligned}
\]
where $\mathcal H$ denotes the Hilbert transform (${\rm p.v.}=$principal value, $\ast$ is the convolution, and $(\cdot)^{\vee}$ is the inverse Fourier Transform)
\[
\begin{aligned}
%\begin{split}
\mathcal H f(x) :=&~ \frac{1}{\pi} {\rm p.v.}\Big(\frac{1}{x}\ast f\Big)(x)
\\
:=& ~ \frac{1}{\pi}\lim_{\epsilon\downarrow 0}\int\limits_{|y|\ge \epsilon} \frac{f(x-y)}{y}\,dy = (-i\,\sgn(\xi) \widehat{f}(\xi))^{\vee}(x).
%\end{split}
\end{aligned}
\]
This model is globally well-posed in the energy space $H^{1/2}$, but also in $L^2$ (Ionescu-Kenig \cite{IK}). Once again, for a complete description of this model, the reader can consult Linares-Ponce \cite{LP}.

\medskip

Now we relax the hypothesis required in the KdV case. Assume that:  
\ben
\item $u=u(t,x)\in C(\R:H^{1/2}(\R)) \cap L^{\infty}_{loc}(\R:L^1(\R))$ solves (BO).

\item  There are $a\in [0,1/2)$ and $c_0>0$ such that for all $T>0$,
\[
\sup_{t\in[0,T]} \int_{\R} |u(t,x)|dx\leq c_0\,\langle T\rangle ^a, \qquad\;\;\;\;\;\; \langle T\rangle:= (1+T^2)^{1/2}.
\]
\een
Therefore, the solution may only have a controlled growth of the $L^1_x$ norm measured in terms of the time variable $T$.

\medskip

\noindent
{\it
{\bf Theorem (G. Ponce, M. \cite{MuPo2}).} Under the previous hypotheses, one has
\[
\liminf_{t\uparrow \infty}\; \int_{\R} \frac{ (u^2+(D^{1/2}u)^2)(t,x)}{1+\left(\frac{x}{\lambda(t)}\right)^2}\,dx= 0,
\]
with
\[
\lambda(t)=\frac{c\,t^{1-a}}{\log t },
\]
for any fixed $\,c>0$.
}

\medskip

The previous result states that at least for a sequence of times, the energy norm in BO is converging to zero in the growing interval $|x| \lesssim t^{(1-a)-}.$

\medskip

Although the proof of this result is similar to the KdV case, it contains several difficulties related to the nonlocal character of the equation. Some previous estimates by Kenig and Martel \cite{MR2590690} are key to prove the final estimates. Some fine commutator estimates are needed, and the liminf is obtained since we need a particular modification of the argument in \cite{MuPo1} that allows us to control the energy norm as time tends to infinity. 

\medskip

\noindent
{\bf More applications:} The method of proof above is quite general and it can be applied to several other dispersive models of importance: BBM equations ({\color{black} Kwak-M. \cite{KM1}}), $abcd$ type models \cite{KMPP,KwM} (small data case only, large data seems quite complicated), Camassa-Holm type equations (Alejo-Cortez-Kwak-M. \cite{ACKM}), KdV on the half line \cite{CM}, the quite involved Zakharov-Rubenchik / Benney-Roskes system (Mar\'ia Mart\'inez and Jos\'e Palacios \cite{MaPa}), and the ILW model ({\color{black} M.-Ponce-Saut} \cite{MPS}).

\section{Coming back to KdV}

The decay result for the BO model has been subsequently improved in the following sense. Let us discuss the recent result by {\color{black} Mendez-M.-Poblete-Pozo} \cite{MMPP}. Let
\[
\Omega(t):=\left\{ x\in\mathbb R ~ :~ |x \pm t^n | \lesssim t^b \right\},  \quad {\color{black} b< \frac{2}{3}}, \quad {\color{black}  0\leq n\leq 1-\frac{b}2}.
\]
For this set, we have 

\medskip

{\it
\noindent
{\color{black}
{\bf Theorem.} Suppose that $u_0\in L^2(\mathbb R)$. Let $u=u(t,x)$ be the corresponding solution to KdV. Then
	\[
	\liminf_{t\rightarrow \infty}\int_{\Omega(t)}u^{2}(t,x)\mathrm{d}x=0. \quad \hbox{(Decay)}
	\]
}
}
This new result improves the previous ones in the following sense: first, it does not require data in $L^1$, only in $L^2$. Second, it expands the region of space where the decay may hold, up to $t^{2/3-}$, in the noncentered case. Third, if one considers decay along rays $|x|\sim t^n$, and $n$ is bigger than $2/3$, then the window width $b$ must be decreased according to the relation $b\leq 2(1-n)$. 

\medskip

\noindent
{\bf Some additional remarks.}
\ben
\item This result is still valid for nonintegrable perturbations of KdV of the form $u^2 + o(u^2)$ and data in $H^1$.
\item The proof also works for quartic gKdV ($u^4$) but we need data in $H^1$. One has
\[
	\liminf_{t\rightarrow \infty}\int_{\widetilde\Omega(t)}u^{4}(t,x)\mathrm{d}x=0, \quad 
\]
with $\widetilde\Omega(t)$ given by
\[
\widetilde\Omega(t):=\left\{ x\in\mathbb R ~ :~ |x \pm t^n | \lesssim  t^b \right\},  \quad {\color{black} b< \frac{4}{7}}, \quad {\color{black}  0\leq n\leq 1-\frac{b}2}.
\]
\item Sharp result in view of Airy decay estimates established by {\color{black} Ifrim-Koch-Tataru} \cite{IKT}.  Indeed, if one assumes ``linear behavior'' for the solution at large time, then (Decay) is sharp.
\een

Obtaining the remaining {\color{black} $\limsup$ property is an open question}, even in the small data case. For further literature in this direction, see the recent result by Duyckaerts, Kenig, Jia and Merle 2017 \cite{DKJM} involving convergence along a sequence of times in the soliton resolution conjecture for the focusing, energy critical wave equation.

\medskip

\section{More applications: the Zakharov-Kuznetsov model}

The previous technique is stable under perturbations: one can prove asymptotic decay for solutions to several dispersive models in more than one dimension.

\medskip

As an application, consider the {\bf Zakharov-Kuznetsov model (ZK)}  
\[
{\color{black}
\hbox{(ZK)} \qquad \partial_t u+\partial_x\Delta u+u\partial_x u=0,
}
\]
where $u=u(t,{\bf x})\in\mathbb{R}$, $t\in\mathbb{R}$ and ${\bf x}\in\mathbb{R}^d$, with $d=2,3$.

\medskip

Concerning this model, one has the following basic facts:
\ben
\item It was originally proposed by Zakharov and Kuznetsov \cite{ZK} (1974) in 3D (see Lannes-Linares-Saut \cite{Lannes:Linares:Saut:2013} for a rigorous derivation). 
\medskip
\item It was derived as an asymptotic model of wave propagation in a magnetized plasma \cite{ZK}.
\medskip
\item It is also a natural multi-dimensional generalization of the Korteweg-de Vries (KdV) equation:
\een
\[
\partial_t u + \partial_x (\partial_x^2 u + u^2)=0, \quad u=u(t,x)\in\R, \quad t,x\in \R.
\]

Formally, it has mass and energy as conserved quantities at the $H^1$ level:
\[
\begin{aligned}
& \int_{\mathbb R^d} u^2(t,{\bf x} )d{\bf x}=\hbox{const.}, \\
& \frac12 \int_{\mathbb R^d} |\nabla u|^2(t,{\bf x} )d{\bf x} - \frac13 \int_{\mathbb R^d} u^3(t,{\bf x} )d{\bf x}=\hbox{const.}
\end{aligned}
\]
Let us review the Cauchy problem for ZK. First of all,

\medskip

\noindent
{\bf $\bullet$ ZK 2D is globally well-posed (LWP) in $L^2$}: See the works by Faminskii \cite{Faminskii:1995}, Linares and Pastor \cite{Linares:Pastor:2009}, Molinet and Pilod \cite{Molinet:Pilod:2015}, Gr\"unrock and Herr \cite{Grunrock:Herr:2015}, and very recently by  Kinoshita \cite{Kinoshita:Arxiv:2019}. Up to know, the best local well-posedness result is in $H^s$, $s>-\frac14$.

\medskip

\noindent
{\bf $\bullet$ ZK 3D is globally well-posed in $H^1$}: Here one has Linares and Saut \cite{MR2486590}, Ribaud-Vento \cite{Ribaud:Vento:2012}, Molinet and Pilod \cite{Molinet:Pilod:2015}, Herr and Kinoshita \cite{Herr:Kinoshita:Arxiv:2020}; the LWP holds for $s>-\frac{1}{2}$.

\medskip

Other results: see the recent uniqueness results vs. spatial decay (Cossetti-Fanelli-Linares \cite{MR3946612}), and the propagation of regularity along regions of space (Linares-Ponce \cite{MR3842873}, Mendez \cite{Mendez}).

\medskip

\noindent
{\bf Solitons.} Similar to KdV, (ZK) possesses soliton solutions of the form
\[
{\color{black} u({\bf x},t)= Q_c(x-ct,x'), \quad c>0, \quad x'\in \mathbb R^{d-1}, \quad d\leq 4.}
\]
Here $Q_c= c Q(\sqrt{c} {\bf x})$ and $Q$ is the $H^1(\mathbb R^d)$ radial solution of the elliptic PDE
\[
{\color{black} \Delta Q -Q +Q^2=0, \quad Q>0.}
\]
Unlike KdV, no explicit formula is known for ZK solitons. However, the following important results are valid:

\ben
\item {Anne de Bouard \cite{MR1378834} (1996):} subcritical ZK solitons are orbitally stable in $H^1$, and supercritical ones are unstable. 
\medskip
\item {C\^ote-M.-Pilod-Simpson \cite{Cote:Munoz:Pilod:Simpson:2016} (2016):} Asymptotic stability (AS) of 2D solitons in the energy space $H^1$.
\medskip
\item Farah-Holmer-Roudenko-Yang (2020): AS in the more involved 3D case \cite{Farah1}. 
\een

Both last works are nontrivial extensions of the foundational works by Martel and Merle concerning the one dimensional KdV case. 

\medskip

\noindent
{More results about solitons:}

\ben
\item Well-decoupled multi-solitons were proved $H^1$ stable in 2D, see \cite{Cote:Munoz:Pilod:Simpson:2016}. 
\medskip
\item  The modified ZK equation (cubic nonlinearity) has finite or infinite time blow up solutions (Farah-Holmer-Roudenko-Yang, \cite{Farah2}). See also Merle \cite{Merle} for the proof in the gKdV case.
\medskip
\item  Construction of multi-solitons in 2D and 3D (Valet, \cite{Valet}), following Martel \cite{Martel}. 
\een

No scattering seems to be available for 2D and 3D, not even in the small data case, except if the nonlinearity is big enough (Farah-Linares-Pastor \cite{MR2950463}). Indeed, ZK in 2D is {\bf scattering critical} ($u u_x \sim \frac1t u$), and subcritical in 3D. KdV in 1D es scattering supercritical. 

\medskip

Our main goal here is to study the asymptotic behavior of global ZK solutions, under minimal regularity assumptions (essentially, data only in $L^2(\mathbb R^d)$ or $H^1(\mathbb R^d)$). 
In particular, we shall describe the dynamics in local regions of space where solitons are absent, i.e., purely ``linear'' behavior.

\medskip

\begin{figure}[h!]
\begin{center}
\begin{tikzpicture}[scale=0.8]
%\filldraw[thick, color=lightgray!30] (-1,1.5)--(5.2,1.5) -- (5.2,5) --(-1,5) -- (-1,1.5);
%\filldraw[thick, color=lightgray!10] (-1,-1)--(-1,4) -- (4,4) --(4,-1) -- (-1,-1);
%\filldraw[thick, color=lightgray!30] (0,4.7)--(0,2) -- (4/3,1) --(4,3) -- (4,4.7) -- (0,4.7);
%\filldraw[thick, color=lightgray!70] (0,4.7)--(0,2) -- (0.6,1.7) --(1,1.6) --(4/3,1.6) --(1.97,1.6) --(4,3) -- (4,4.7) -- (0,4.7);
%%\draw[thick, color=black] (1.25,4) -- (3.2,1.5);
\filldraw[thick, color=lightgray!30] (0,0) -- (4,0) -- (4,4) -- (0,4) -- (0,0);
\draw[thick, dashed] (0,0) -- (4,0) -- (4,4) -- (0,4) -- (0,0);
%\draw[thick] (4,3) -- (4,5);
\draw[thick,dashed] (5,-1) -- (5,5);
%\draw[thick,dashed] (0,2)--(8/3,0);
%\draw[thick] (0,2) -- (0,5);
%\draw[thick,dashed] (3.2,-1)--(3.2,5.3);
%\draw[thick,dashed] (-1,1.5)--(5.2,1.5);
%\draw[thick,dashed] (-1,5)--(5.2,5);
\draw[->] (-1,2) -- (6,2) node[below] {$x$};
\draw[->] (2,-1) -- (2,5) node[right] {$y$};
%\node at (1.9,-0.7){$ \mathcal{B}_2(b)$};
%\node at (1.2,2.2){$- \frac{1+4c}{2(1-c)}$};
%\node at (0,0){$\bullet$};
%\node at (0,4){$\bullet$};
%\node at (4,0){$\bullet$};
\node at (5.7,2.5){$x\sim t$};
\node at (2.5,3.7){$t^{br}$};
\node at (3.7,2.5){$t^b$};
\node at (1,1){$\Omega(t)$};
%\node at (4/3,0){$\bullet$};
%\node at (8/3,-0.4){$\frac23$};
%\node at (4/3,-0.4){$\frac13$};
%\node at (8/3,0){$\bullet$};
%\node at (0,2){$\bullet$};
%\node at (-0.3,2){$\frac 13$};
%\node at (0,1){$\bullet$};
%\node at (-0.3,1){$\frac16$};
%%\node at (0,8/3){$\bullet$};
%\node at (0,3){$\bullet$};
%\node at (4,1.6){$\bullet$};
%\node at (4.3,1.6){$\frac29$};
%\node at (-0.3,3){$\frac12$};
%\node at (-1,0.5){$x$};
%\node at (1,0.5){$y$};
%\draw[thick,dashed] (2,4) arc (180:360:0.6);
%\draw[thick,dashed] (4,2) arc (90:270:0.6);
%\draw[thick,dashed] (2.5,2) arc (0:360:0.6);
%\node at (4.7,4.6){$-b+\frac12$};
%\node at (4.7,3.3){$ \mathcal{B}_3(b)$};
%\node at (-3.5,2.5){$V$};
%\node at (3,2.7){$W$};
%\draw (2,0.5) arc (90:145:0.5);
%\node at (2.7,-1.9){$x_1=\beta t$};
%\node at (4,-1.2){$x_1+(\tan \theta) \, x_2=\beta t $};
%\node at (0.3,0.15){$\theta$};
%\draw (0.5,0) arc (0:55:0.5);
%\node at (0.1,1.1){$\theta$};
%\draw (0,0.9) arc (270:325:0.4);
%\node at (-0.7,3.5){$\mathcal{B}_0$};
\end{tikzpicture}
\caption{ Here {\color{black} $\frac13<r<3$, $0\leq b<\frac{2}{3+r}$ and $0\leq br<\frac{2r}{3+r}$}. Here $x\sim t$ represents the soliton region.}\label{fig:1}
\end{center}
\end{figure}
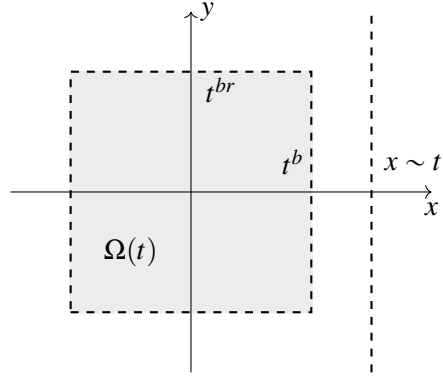

{\it
\noindent
{\bf Theorem (Mendez-M.-Poblete-Pozo 2020 \cite{MMPP}).}
	Suppose that $u_0\in L^2(\mathbb R^2)$ and let $u=u(x,y,t)$ be the bounded in time solution  to 2D (ZK). Then
	\[
	\liminf_{t\rightarrow \infty}\int_{\Omega(t)}u^{2}(x,y,t)\mathrm{d}x\,\mathrm{d}y=0,
	\]
	with 
	\[
	\Omega(t)=\left\{ (x,y)\in\R^2 ~ :~|x| < t^b, \quad |y|<t^{br}, \quad \frac13<r<3, \quad 0\leq b<\frac{2}{3+r}\right\}.
	\]
	See Fig. \ref{fig:1}.
}	
%	Moreover, there exist constant $ C_0>0$ and an increasing sequence of times $t_n\to +\infty$ such that
%	\[
%	\int_{\Omega(t_n)}u^{2}(x,y,t_n)\mathrm{d}x\,\mathrm{d}y\leq \frac{C_0}{\ln^{2} (t_n) }.
%	\]
\medskip

%\begin{figure}
%  \begin{minipage}[c]{0.67\textwidth}
%    \includegraphics[width=\textwidth]{2011-03-03}
%  \end{minipage}\hfill
%  \begin{minipage}[c]{0.3\textwidth}
%    \caption{
%       ??????????? ???????? ? ???????????? ?????????
%       ?????????? ???????? ? ??????????? ?? ????????? ????.
%    } \label{fig:03-03}
%  \end{minipage}
%\end{figure}

{\bf Some comments.}
\ben
\item The previous result also valid in 3D and if data is in $H^1$. One can also show that the $L^2$ local gradient decays to zero.
\medskip
\item One can also show strong $L^2$ decay in regions of the form $|x|\gg t$, see {\color{black} M.-Ponce-Saut} \cite{MPS} for the case of the ILW equation.
\medskip
\item See another application to the {KdV-Schr\"odinger} model by Linares-Mendez \cite{LM}.
\een
As we can see, the virial technique is stable enough to consider 2 and 3 dimensional cases. The remaining regions in ZK are of independent interest, in view that there is a soliton region $x\sim t$.

\section{Coming back again to gKdV}

The following result is related to the limsup problem described above. A recent result shows complete decay to zero in extreme regions of space. Of particular interest is the left region in the gKdV case, which was completely out of reach.

\medskip

\noindent
{\it
{\bf Theorem (R. Freire, G. Ponce, F. Linares, M.).} Any global solution to subcritical gKdV in $L^\infty_t H^1_x$ satisfies
\[
\lim_{t\to +\infty} \| u(t)\|_{L^2(\Omega(t))} =0,
\]
with
\[
\Omega(t):= \left\{ x\in \R ~ :~ x\lesssim -t\log^{1+\epsilon} t \quad \hbox{and}\quad x> C_0t \right\},
\]
with $C_0$ depending on the size of the initial data and any $\epsilon>0$. A similar result also holds for the BO equation.
}

\medskip

Note that this result also holds for the cubic power $p=3$, which has breathers. In some sense, the previous result is sharp, since breathers \eqref{breather} do not decay and can move along lines $x=-vt$, for any $v>0$.

\medskip

With this result in mind, one has a satisfactory answer in the nonintegrable case and in far regions of space. The central region is well-known only in terms of liminf, and the soliton region in full generality is still completely out of reach.

\medskip

The proof of these results is done by using the mass of the solution cut into particularly well-chosen pieces. See \cite{FLMP} for further details.

\section{The KP models}

The {\bf Kadomtsev-Petviashvili (KP)} equations in $\R^2$ are:
{\color{black}
\[
\partial_t u+\partial_{x}^{3}u +u\partial_x u +\kappa\partial_{x}^{-1}\partial_{y}^{2}u=0,
\]
}
where $u=u(t,x,y)\in\mathbb{R}$, $t\in\mathbb{R}$, $(x,y)\in\mathbb{R}^2$ and  $\kappa\in \{-1,1\}.$

\medskip

The nonlocal operator $\partial_{x}^{-1}f$ is formally defined in the literature as
\begin{equation*}
(\partial_{x}^{-1}f)(x,y):=\int_{-\infty}^{x}f(s,y)\,\mathrm{d}s.
\end{equation*}
The KP equations were first introduced by Kadomtsev and Petviashvili in 1970 \cite{KadomtsevPetviashvili1970} for modeling {long and weakly nonlinear waves} propagating essentially along the $x$ direction, {\color{black} with a small dependence in the $y$ variable}. A rigorous derivation of both models from the symmetric $abcd$ Boussinesq system was obtained by {Lannes and Lannes-Saut \cite{Lannes2,LS}.} The nonlocal term makes KP models {hard from the mathematical point of view}, and the understanding of the dynamics is far from satisfactory.

\medskip

The convention here is that $\kappa=-1$ is KP-I (strong surface tension) and $\kappa=1$ is KP-II (weak surface tension). 

\medskip

Despite their apparent similarity, KP-I and KP-II  differ significantly with respect to their underlying mathematical structure and the behavior of their solutions. For instance, from the point of view of {well-posedness} theory KP-II is much better understood than KP-I.

\medskip

\noindent
{\bf KP II}. Bourgain \cite{Bourgain1993} showed that KP-II is globally well-posed (GWP) in $L^2(\mathbb{R}^2)$ (see also Ukai \cite{Ukai} and I\'orio-Nunes \cite{IN} for early results). Improvements of Bourgain results by Takaoka-Tzvetkov \cite{Takaoka-Tzvetkov-2001}, Isaza-Mejia \cite{Isaza-Mejia-2001} and Hadac-Herr-Koch \cite{Hadac-Herr-Koch-2009}.
%\een

\medskip

Bourgain proved GWP via the contraction principle in $X^{s,b}$ spaces and the conservation of the $L^2$-norm
\[
{\color{black} M[u](t):= \frac12\int_{\mathbb R^2} u^2(t,x,y)\,\mathrm{d}x\,\mathrm{d}y = \hbox{const.}}
\]
Note that this conservation holds for both KP-I and KP-II.

\medskip

\noindent
{\bf KP I}. The KP-I global theory in the energy space took years to be solved. It was known from Molinet-Saut-Tzvetkov \cite{Molinet-Saut-Tzvetkov-2002_b} that KP-I {behaves badly with respect to perturbative methods}.

\ben
\item GWP by Molinet, Saut and Tzvetkov \cite{Molinet-Saut-Tzvetkov-2002}, C.~E. Kenig  \cite{Kenig-2004}, and Ionescu, Kenig and Tataru \cite{Ionescu-Kenig-Tataru-2008}. See also \cite{Ionescu-Kenig-2007}. 
\medskip
\item The last result is obtained in the natural {\bf energy space} of the equation 
\[
{\color{black}
E^1(\mathbb R^2) := \left\{ u \in L^2 (\mathbb R^2) ~ : ~ 
\|u \|_{L^2}+\|\partial_xu\|_{L^2} +\|\partial_x^{-1}\partial_yu \|_{L^2} <+\infty \right\},
}\]
in the sense that the flow map extends continuously from suitable smooth data into $E^1(\mathbb R^2)$.
\een 
The space $E^1(\mathbb R^2)$ arises naturally from the {\bf conservation of energy} ($\kappa=1$ KP-II, $\kappa=-1$ KP-I)
\[
{\color{black}
E[u](t):= \int_{\mathbb R^2} \left(\frac12 (\partial_x u)^2 -\frac12\kappa  (\partial_x^{-1}\partial_y u)^2 -\frac13 u^3 \right)(t,x,y)\,\mathrm{d}x\,\mathrm{d}y = \hbox{const.}
}\]
Another important conserved quantity in KP models is the momentum
\[
{\color{black}
P[u](t):= \frac12\int_{\mathbb R^2} \left( u\partial_x^{-1}\partial_y u \right) (t,x,y)\,\mathrm{d}x\,\mathrm{d}y = \hbox{const.}
}\]

%\begin{frame}
%
%{\footnotesize
%
\noindent
{\bf Lumps}. {\color{black} KP-I has lump solutions}, namely solutions of the form
\[
u(t,x,y)=Q_c(x-ct,y), \quad c>0.
\]
The function $Q_c$ is given as $Q_c(x,y):=cQ(\sqrt{c} x,c y)$, where $Q$ is the fixed profile
\[
{\color{black}
Q(x,y)=12\partial_x^2 \log(3+x^2+y^2)=\frac{24(3-x^2+y^2)}{(x^2+y^2+3)^2}.
}\]
This profile satisfies the elliptic nonlocal PDE on $\mathbb R^2$
\[
{\color{black}
\partial_x^2 Q -Q + \frac12 Q^2 -\partial_x^{-2}\partial_y^2 Q=0, \quad Q\in E^1(\mathbb R^2),
}\]
in the sense of distributions. Lumps were first found by Satsuma and Ablowitz via an intricate limiting process of complex-valued algebraic solutions to KP-I \cite{SatAblo}. 

\medskip

\noindent
{\bf Lumps vs. ground states}. De Bouard and Saut \cite{dBS,dBS1,dBS2} described, via PDE techniques, qualitative properties of KP-I lumps solutions which are also ground states (not necessarily equal to $Q$). 

\medskip

Ground states decay as $1/(x^2+y^2)$ as $|(x,y)|$ tends to infinity. However, whether or not lump solutions $Q$ are ground states {\color{black} is still an unknown open problem}.

\medskip

Liu and Wei \cite{LiuWei}, by using  linear B\"acklund transformation techniques, proved the {\bf orbital stability of the lump $Q$ in the space $E^1(\mathbb R^2)$}, hinting that $Q$ it is probably the unique (modulo translations) ground state. 
%
%}
%
%\end{frame}

\medskip

%
%\begin{frame}
%
%{\footnotesize
%
%
Let us review the invariances present in the KP equations:
\begin{enumerate}
\item Shifts: 
\[
u(t,x,y) \mapsto u(t+t_0, x+x_0,y+y_0).
\]
\item Scaling: if $c>0$,
\[
u(t,x,y) \mapsto c u\left( c^{3/2} t , c^{1/2} x, c y\right).
\]
\item Galilean invariance: for any $\beta\in\mathbb R$, if $u(x,y,t)$ is solution to KP, then
\[
\begin{cases}
u(t,x,y) \mapsto u\left(t,x -\beta^2 t -\beta (y-2\beta t) , y-2\beta t \right), & \hbox{KP-I}\\
u(t,x,y) \mapsto u\left(t,x + \beta^2 t + \beta (y-2\beta t) , y-2\beta t\right), & \hbox{KP-II},
\end{cases}
\] 
define new solutions to KP.
\end{enumerate}

Using this, one can construct a moving lump solution, of arbitrary size and speed: for any $\alpha,\beta\in\mathbb R$,
\[
{\color{black}
Q_{c,\beta}(t,x,y):= Q_c \left(x-ct -\beta^2 t -\beta (y-2\beta t) ,y-2\beta t \right)
}
\]
is a moving lump solution of KP-I with speed $(c+\beta^2, 2\beta)\in\mathbb R^2$. One can show that small lumps (in the energy space) may have arbitrarily large speeds ($\beta\gg1$ and $c\ll 1$ such that $\beta^2 c^{1/2}\ll 1$). %Precisely, this property of lumps ensures that \emph{no monotonicity on the right of the solution} may hold in the KP-I case, in the sense of Martel and Merle. % (note that this is not the case in KP-II, see \eqref{zerozero} and de Bouard-Martel \cite{dBM}). %A general solution may contain a small fast lump solution that will appear after some time on the right of the main part of the solution. This simply implies that an arbitrary portion of mass on the right of the plane cannot be almost preserved in time. In that sense, the proof of Theorem \ref{thm_KPI} overcomes this difficulty by adding one more degree of regularity, and using the momentum as a trigger of decay, instead of considering only the mass as in standard monotonicity properties.

%
%}
%
%\end{frame}

\medskip

%
%\begin{frame}
%
%{\footnotesize
%
In the KP-II case:

\ben
\item de Bouard-Saut \cite{dBS1}, and de Bouard-Martel \cite{dBM}: KP-II has no lump structures as KP-I has. 
\medskip
\item Any KdV soliton becomes an (infinite energy) line-soliton solution of KP. This structure is {\bf stable} in the KP-II case, as proved by Mizumachi and Tzvetkov \cite{MT}, and {\bf asymptotically stable} in a series of deep works by Mizumachi \cite{Mizu1,Mizu2}. 
\medskip
\item  Multi-line-soliton structures are known to exist via Inverse Scattering methods, and nothing is known about their stability in rigorous terms \cite{KS}.  
\medskip
\item In the KP-I case, the KdV soliton is transversally unstable (Rousset and Tzvetkov \cite{RT1,RT2}). 
\medskip
\item Molinet, Saut and Tzvetkov \cite{Molinet-Saut-Tzvetkov-2011} proved global well-posedness of KP-II along the KdV line soliton in $L^2(\mathbb{R}\times \mathbb{T})$  and $L^2(\mathbb R^2)$ (see also Koch and Tzvetkov \cite{Koch-Tzvetkov-2008}). 
\medskip
\item Izasa-Linares-Ponce \cite{ILP_2} showed propagation of regularity for this model. The KP-I case seems unsolved, as far as we understand.
\een
%}
%
%\end{frame}
%
%

Once a suitable well-posedness theory is available, one may wonder about the {\bf long time behavior of globally defined solutions}. This is a difficult open question not yet solved either by Inverse Scattering Transform (IST) or PDE methods (see \cite{W,Sung} for partial results in this direction, in the case of small data). Moreover, the answer may be strongly dependent on the choice of model, KP-I or KP-II.

\medskip

Except by some particular cases ({lumps and line solitons}, their orbital and asymptotic stability, and the case of suitable {small data scattering solutions}), no rigorous results about large data behavior in KP models, starting from Cauchy data, are available, as explained by Klein-Saut \cite{KS}.

\begin{figure}[h!]
  \centering
  \includegraphics[width=5cm]{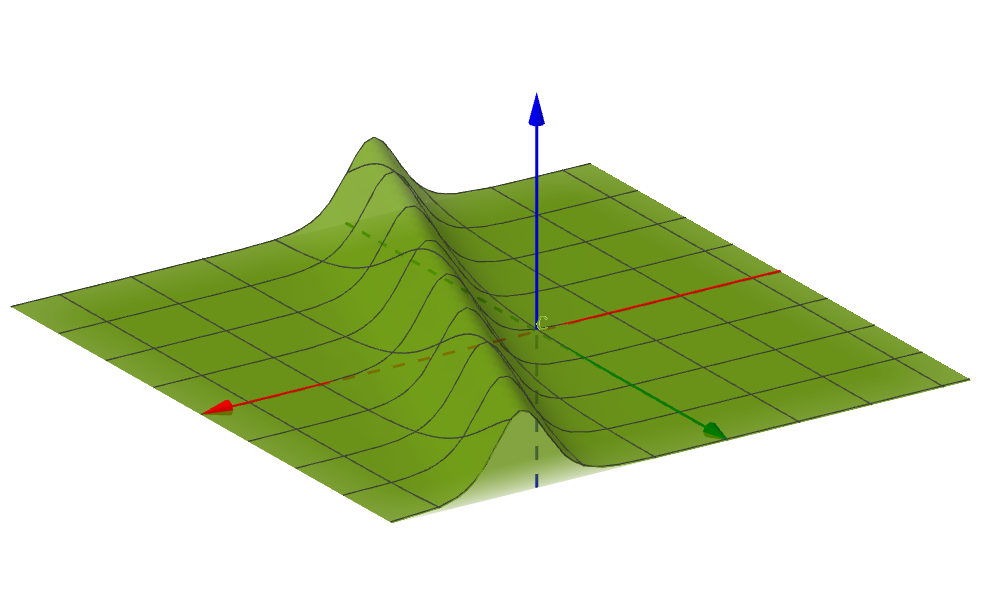}
  \includegraphics[width=5cm]{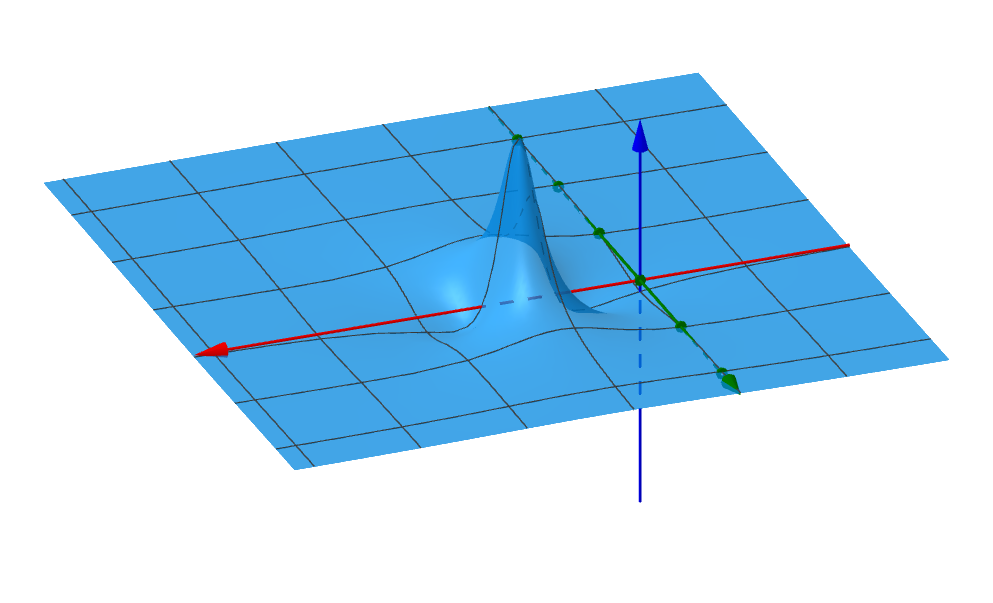}
  \caption{Line Soliton of KP (left) and KP-I lump (right).}
\end{figure}

\medskip

\noindent
{\bf New results.} Let $t\gg 1$ and $\Omega_1(t)$ denote the following rectangular box 
\[
{\color{black}
\Omega_1(t)= \left\{ (x,y)\in\mathbb R^2 ~: ~  |x-\ell_1 t^{m_1}| \leq t^{b}, \; |y-\ell_2 t^{m_2}|\leq t^{br}\right\},
}\]
with $\ell_1,\ell_2\in\mathbb R$,
{\color{black}
\[
\begin{aligned}
&\frac53<r<3, \quad 0<b <\frac2{3+r}, \\
& 0\leq m_1<1-\frac{b}2(r+1), \quad \hbox{and}\quad  0\leq m_2 < 1-\frac{b}2(3-r).
\end{aligned}
\]
}
Let also $\sigma_1,\sigma_2\geq 0$, not both equal zero. We define $\Omega_{2,1}(t)$ as 
{\color{black}
\[
\Omega_{2,1}(t):= \left\{ (x,y)\in\mathbb R^2 ~: ~   \sigma_1  \left| x \right| +\sigma_2  \left| y \right| \geq t \log^{1+\epsilon }t \right\},
\]
}for any fixed $\epsilon>0.$ 

\medskip

We also define, for $\beta>0$ and $\sigma_3\in\mathbb R$,
{\color{black}
\[
\Omega_{2,2}(t):= \left\{ (x,y)\in\mathbb R^2 ~: ~   x +\sigma_3 y \geq \beta  t  \right\}.
\]
}
Under these definitions, we proved:
%Let also $\Omega_2(t)$ be the union of cylinders 
%\[
%{\color{black}
%\Omega_2(t):= \left\{ (x,y)\in\mathbb R^2 ~: ~  |x|\sim t^{p}\log^{1+\epsilon }t \quad  \hbox{ or } \quad  |y|\sim t^{p}\log^{1+\epsilon }t \right\},
%}\]
%for any fixed $p\geq 1$ and $\epsilon>0.$ Here $a\sim b$ means $c_0 a\leq b\leq C_0 a$, for some fixed $c_0,C_0>0.$ %See Fig. \ref{fig:1} for details.

\medskip

\begin{figure}[h!]
\begin{center}
\begin{tikzpicture}[scale=0.7]
%\filldraw[thick, color=lightgray!30] (-1,1.5)--(5.2,1.5) -- (5.2,5) --(-1,5) -- (-1,1.5);
%\filldraw[thick, color=lightgray!10] (-1,-1)--(-1,4) -- (4,4) --(4,-1) -- (-1,-1);
%\filldraw[thick, color=lightgray!30] (0,4.7)--(0,2) -- (4/3,1) --(4,3) -- (4,4.7) -- (0,4.7);
%\filldraw[thick, color=lightgray!70] (0,4.7)--(0,2) -- (0.6,1.7) --(1,1.6) --(4/3,1.6) --(1.97,1.6) --(4,3) -- (4,4.7) -- (0,4.7);
%%\draw[thick, color=black] (1.25,4) -- (3.2,1.5);
\filldraw[thick, color=lightgray!30] (0,-1) -- (4,-1) -- (4,5) -- (0,5) -- (0,-1);
\draw[thick, dashed] (0,-1) -- (4,-1) -- (4,5) -- (0,5) -- (0,-1);
%\draw[thick] (4,3) -- (4,5);
\draw[thick,dashed] (-1,5.5)--(5,5.5);
\draw[thick,dashed] (-1,-1.5)--(-1,5.5);
\draw[thick,dashed] (-1,-1.5)--(5,-1.5);
%\draw[thick] (0,2) -- (0,5);
%\draw[thick,dashed] (3.2,-1)--(3.2,5.3);
%\draw[thick,dashed] (-1,1.5)--(5.2,1.5);
%\draw[thick,dashed] (-1,5)--(5.2,5);
%(4 sigueintes del interior para gradientes)
%\filldraw[thick, color=lightgray!60] (0.1,0.2) -- (3.9,0.2) -- (3.9,3.8) -- (0.1,3.8) -- (0.1,0.2);
%\draw[thick, dashed] (0.1,0.2) -- (3.9,0.2) -- (3.9,3.8) -- (0.1,3.8) -- (0.1,0.2);
%\filldraw[thick, color=lightgray!90] (0.3,0.4) -- (3.7,0.4) -- (3.7,3.6) -- (0.3,3.6) -- (0.3,0.4);
%\draw[thick, dashed] (0.3,0.4) -- (3.7,0.4) -- (3.7,3.6) -- (0.3,3.6) -- (0.3,0.4);
\filldraw[thick, color=lightgray!30] (6,-3) -- (6,7) -- (7,7) -- (7,-3) -- (6,-3);
\filldraw[thick, color=lightgray!30] (-3,-3) -- (-3,7) -- (-2,7) -- (-2,-3) -- (-3,-3);
\filldraw[thick, color=lightgray!30] (7,6) -- (7,7) -- (-3,7) -- (-3,6) -- (7,6);
\filldraw[thick, color=lightgray!30] (-3,-2) -- (-3,-3) -- (7,-3) -- (7,-2) -- (-3,-2);
%\draw[thick, dashed] (0.3,0.4) -- (3.7,0.4) -- (3.7,3.6) -- (0.3,3.6) -- (0.3,0.4);
%\filldraw[thick, color=lightgray!30] (0.3,0.4) -- (3.7,0.4) -- (3.7,3.6) -- (0.3,3.6) -- (0.3,0.4);
%\draw[thick, dashed] (0.3,0.4) -- (3.7,0.4) -- (3.7,3.6) -- (0.3,3.6) -- (0.3,0.4);
\draw[->] (-3,2) -- (7.5,2) node[below] {$x$};
\draw[->] (2,-3.5) -- (2,7.5) node[right] {$y$};
%\node at (1.9,-0.7){$ \mathcal{B}_2(b)$};
%\node at (1.2,2.2){$- \frac{1+4c}{2(1-c)}$};
\node at (6.5,2.3){$x\gg t$};
\node at (2.7,6.5){$y\gg t$};
%\node at (0,4){$\bullet$};
%\node at (4,0){$\bullet$};
\node at (5.2,2.5){$x\sim t$};
\node at (1,5.8){$y\sim t$};
\node at (2.4,4.7){$t^{br}$};
\node at (4.3,2.5){$t^b$};
\node at (0.6,-0.5){$\Omega_1(t)$};
\node at (6.5,3){$\Omega_{2,1}(t)$};
\node at (6,0){$\Omega_{2,2}(t)$};
%\node at (3.4,-0.2){$\widetilde\Omega_1(t)$};
%\node at (0.9,3.1){$\widetilde\Omega_{2}(t)$};
\draw[thick,dashed] (5,-3) -- (5,7);
\end{tikzpicture}
\end{center}
%\caption{\small Schematic figure depicting the sets $\Omega_1(t)$, $\Omega_{2,1}(t)$, $\Omega_{2,2}(t)$, $\widetilde \Omega_1(t)$ and $\widetilde \Omega_2(t)$, \emph{in the centered case} $l_1=l_2=0$, as defined in \eqref{Omega}, \eqref{Omega0}, \eqref{Omega0b}, \eqref{Omega1} and \eqref{Omega2}, respectively. Recall that $\frac53<r<3$, $0< b<\frac{2}{3+r}$ and $0< br<\frac{2r}{3+r}$. The largest value of $b$ is $\sim\frac{3}{7}$, and $br\sim \frac57$, obtained by $r\sim \frac53$. Most of the white area outside the regions $|x|\sim t$ and $|y|\sim t$ \emph{can be suitably covered} by using shifts, in the non-centered case.}\label{fig:1}
\end{figure}

\noindent
{\it
{\bf Theorem \cite{MMPP1}.}
Every solution $u=u(x,y,t)$ of KP obtained from arbitrary initial data $u_0$ in $L^2(\mathbb R^2)$ in the KP-II case, and $u_0$ in the energy space $E^1(\mathbb R^2)$ for KP-I, satisfies% Then the corresponding solution $u$ satisfies
\[
 \liminf_{t\to \infty} \int_{\Omega_1(t)} u^2(t,x,y)\,\mathrm{d}x\mathrm{d}y=  0.
\]
If now $u\in L^\infty([0,\infty); E^1(\mathbb R^2))$ is a solution to KP,
\[
 \lim_{t\to \infty} \int_{\Omega_{2,1}(t)} u^2(t,x,y)\,\mathrm{d}x\mathrm{d}y= 0,
\]
and for KP-II and $\beta$ sufficiently large (only depending on $\sigma_3$ and the $E^1(\mathbb R^2)$ norm of the initial datum),
\[
 \lim_{t\to \infty} \int_{\Omega_{2,2}(t)} u^2(t,x,y)\,\mathrm{d}x\mathrm{d}y= 0.
\]
%and
%\[
% \lim_{t\to \infty} \int_{\Omega_2(t)} u^2(x,y,t)\,\mathrm{d}x\mathrm{d}y=  0.
%\]
}

\noindent
{\bf Final remarks:}
\ben
\item This result is satisfied by both KP-I and KP-II: both models contain a quadratic nonlinearity, and the proof does not depend on the sign of $\kappa$. 
\medskip
\item In the KP-II setting, more can be said: Kenig and Martel \cite{KM} showed that for any $\beta>0$ and  initial data small in $L^1\cap L^2$,
\[
\lim_{t\to \infty} \int_{x>\beta t} u^2(t,x,y)\,\mathrm{d}x\mathrm{d}y=  0.
\] 
\item The remaining limsup property holds in a particular region if the solution is in $L^\infty_t (L^1_{x,y}\cap L^2_{x,y})$, see Mu\~noz-Ponce \cite{MuPo1} for similar arguments. 
\medskip
\item There is ground to believe that the results in the previous theorem, as well as in the KdV case, are sharp. 
\medskip
\item The previous theorem is stable under perturbations of the nonlinearity $f(u)= u^2 + o_{u\to 0}(u^2)$, provided LWP of the corresponding KP equation is available. 
\een

\end{document}